\documentclass[11pt]{article}
\begin{document}
\begin{center}
{\Large\bf A New Approach to Signed Eulerian Numbers}\\
\vspace{0.5cm}
{\large\bf Shinji Tanimoto}\\
(tanimoto@cc.kochi-wu.ac.jp)\\
\vspace{0.5cm}
Department of Mathematics\\
Kochi Joshi University\\
Kochi 780-8515, 
Japan. \\
\end{center}
\begin{abstract}
The numbers of even and odd permutations with a given ascent number are investigated
by an operator that was introduced in [9]. Their difference is called a signed Eulerian number.        
By means of the operator the recurrence relation for signed  Eulerian numbers can be deduced, 
which was obtained in [1] by an analytic method.
Our approach is straightforward and enables us to deduce other properties including divisibility properties 
by prime powers. 
\end{abstract}
\vspace{1.0cm}
\noindent
{\bf 1.  Introduction} \\
\\
\indent
An {\it ascent} (or {\it descent}) of a permutation $a_1a_2\cdots a_n$ of $[n]=\{1,2, \ldots, n\}$ 
is an adjacent pair such that $a_i < a_{i+1}$ (or $a_i > a_{i+1}$) for some $i$ ($1\le i \le n-1$).  
Let $E(n,k)$ be the set of all permutations of $[n]$ 
with exactly $k$ ascents, where $0 \le k \le n-1$. Its cardinality is the classical
{\it Eulerian number};
\[
    A_{n,k}= |E(n,k)|, 
\]
whose properties and identities can be found in [2-6], for example.\\
\indent
An {\it inversion} of a permutation $A=a_1a_2\cdots a_n$ is a 
pair $(i, j)$ such that $1 \le i < j \le n$ and $a_i > a_j$. Let us denote by ${\rm inv}(A)$     
the number of inversions in a permutation $A$, and by $E_{\rm e}(n,k)$ or $E_{\rm o}(n,k)$ 
the subsets of all permutations in $E(n,k)$ 
that have, respectively, even or odd numbers of inversions. 
The aim of this paper is to investigate their cardinalities;
\[
     B_{n,k} = |E_{\rm e}(n,k)|~~  {\rm and} ~~  C_{n,k} = |E_{\rm o}(n,k)|.
\]   
\indent
Obviously we have $A_{n,k}=B_{n,k}+C_{n,k}$, while
the differences
\[
        D_{n,k}=B_{n,k}- C_{n,k}
\]
were called {\it signed Eulerian numbers} in [1], where descents of permutations
were considered instead of ascents. Therefore, the identities for $D_{n,k}$ presented here correspond to 
those in [1] that are obtained by replacing $k$ with $n-1-k$. \\
\indent
In order to study these numbers, we make use of an 
operator on permutations in $[n]$, which was introduced in [9]. In the subsequent papers [10] and [11], 
it was shown that the operator plays a relevant role in studying Eulerian numbers.
The operator $\sigma$ is defined by adding one to 
all entries of a permutation and by changing $n+1$ into one. 
However, when $n$ appears at either end of a permutation, it is removed
and one is put at the other end. That is, for a permutation $a_1a_2\cdots a_n$ with $a_i = n$ for some $i$ $(2 \le i\le n-1)$, we have
\begin{itemize}
\item[(i)] $\sigma(a_1a_2\cdots a_n) = b_1b_2\cdots b_n$, 
\end{itemize}
\noindent
where $b_i = a_i +1$ for all $i$ $(1 \le i\le n)$ and $n+1$ is replaced by one. 
And, for a permutation $a_1a_2\cdots a_{n-1}$ of $[n-1]$, we have:
\begin{itemize}
\item[(ii)] $\sigma(a_1a_2\cdots a_{n-1}n) = 1b_1b_2\cdots b_{n-1}$; 
\item[(iii)] $\sigma(na_1a_2\cdots a_{n-1}) = b_1b_2\cdots b_{n-1}1$,
\end{itemize}
\noindent
where $b_i = a_i +1$ for all $i$ $(1 \le i\le n-1)$. We denote by $\sigma^{\ell}A$
the repeated $\ell$ applications of $\sigma$ to a permutation $A$.
\\
\indent
It is obvious that the operator preserves the number
of ascents or descents in a permutation, that is, 
$\sigma A \in E(n,k)$ if and only if $A \in E(n,k)$. Let us observe the number of inversions 
of a permutation when $\sigma$ is applied. \\
\indent
When $n$ appears at either end of a permutation $A=a_{1}a_{2}\cdots a_{n}$
as in (ii) or (iii), it is evident that 
\[
{\rm inv}(\sigma A) = {\rm inv}(A).
\]
\indent
Now let us consider the case (i). When $a_{i}=n$ for some $i$ ($2\le i \le n-1$), we get
$\sigma (a_1a_{2}\cdots a_{n}) =b_1b_2 \cdots b_n$, 
where $b_{i}= 1$ is at the position. In this case, $n-i$ inversions 
$(i, i+1), \ldots, (i, n)$ of $A$ vanish and, in turn, 
$i-1$ inversions $(1, i), \ldots, (i-1, i)$ of $\sigma A$ occur. Hence
the difference between the numbers of inversions is 
\begin{eqnarray}
 {\rm inv}(\sigma A) - {\rm inv}(A) = (i-1) - (n-i) = 2i - (n + 1).
\end{eqnarray}
\indent
Therefore, when $n$ is even, each application of the operator changes the parity 
of permutations as long as $n$ remains in the interior of permutations. 
If $n$ is odd, however, the operator $\sigma$ also preserves the parity of all 
permutations of $[n]$.\\
\indent
For convenience sake we denote by $E_{\rm e}^{-}(n,k)$ and $E_{\rm e}^{+}(n,k)$ the sets of permutations
$a_1a_2\cdots a_n$ in $E_{\rm e}(n,k)$ with $a_1 < a_n$ and $a_1 > a_n$, 
respectively. Similarly, $E_{\rm o}^{-}(n,k)$ and $E_{\rm o}^{+}(n,k)$ denote
those in $E_{\rm o}(n,k)$. 
In $E_{\rm e}^{-}(n,k)$ or $E_{\rm o}^{-}(n,k)$ {\it canonical} permutations 
are those of the form $1a_2a_3 \cdots a_n$, and in 
$E_{\rm e}^{+}(n,k)$ 
or $E_{\rm o}^{+}(n,k)$ are those of the form $a_2a_3 \cdots a_n1$, where
$a_2a_3 \cdots a_n$ is a permutation of $\{2, 3, \ldots, n\}$. \\
\indent
In [8] and references therein, even or odd permutations were classified 
by anti-excedance number, not by the ascent number. An anti-excedance in a permutation
$A=a_1a_2\cdots a_n$ means an inequality $i \ge a_i$.
Recurrence relations were also given for the cardinalities of the sets of even and odd 
permutations that are classified by the anti-excedance number. The recurrence 
relations held for all $n$.
The classification of even or odd permutations by the ascent number seems not so         
simple, as will be seen in the following sections. \\
\indent
The signed Eulerian numbers $D_{n,k} = B_{n,k} - C_{n,k}$,
however, have a recurrence relation that holds for all $n$, 
although it has different expressions according to the parity of $n$.
The relation was conjectured in [7] and an analytic proof for it was given in [1]. 
In Section 4 we will derive it from a quite different point of view based on the properties of the operator $\sigma$.  \\
\\
\noindent
{\bf 2.  The Numbers {\boldmath $B_{n,k}$} and {\boldmath $C_{n,k}$}} \\
\\
\indent
The numbers $B_{n,k}$ and $C_{n,k}$ enjoy some symmetry properties 
according to the values of $n$. The permutation $n \cdots 21 \in E(n, 0)$ has 
$n(n-1)/2$ inversions. Hence the values of $B_{n,0}$ and $C_{n,0}$ are given by
\begin{eqnarray*}
     B_{n,0} = 
        \left\{\begin{array}{rl}
                    1, & \mbox{if  $n~ \equiv$ ~0~ or~ 1 ~($\bmod$ 4), } \\
                    0, & \mbox{if  $n~ \equiv$~ 2 ~or ~3 ~($\bmod$ 4),} \\
                    \end{array}  \right.
\end{eqnarray*}
and
\begin{eqnarray*}
     C_{n,0} = 
        \left\{\begin{array}{rl}
                    0, & \mbox{if  $n~ \equiv$~ 0 ~or ~1 ~($\bmod $ 4), } \\
                    1, & \mbox{if  $n~ \equiv$ ~2 ~or ~3 ~($\bmod$ 4).} \\
                    \end{array}  \right.
\end{eqnarray*}
\indent
For a permutation $A = a_1a_2 \cdots a_n$ we define its reflection 
by $A^{\ast}= a_n  \cdots a_2a_1$.
Using reflected permutations and the parity of $n(n-1)/2$, the following symmetry 
between $B_{n,k}$ and $C_{n,k}$ are easily checked. 
\indent
\begin{itemize}
\item[(i)] $n \equiv$ 0  or 1 ($\bmod$ 4). 
In this case, $A \in E_{\rm e}(n,k)$ if and only if
$A^{\ast} \in E_{\rm e}(n,n-k-1)$, and 
$A \in E_{\rm o}(n,k)$ if and only if
$A^{\ast} \in E_{\rm o}(n,n-k-1)$, so we have
\[
     B_{n,k} = B_{n,n-k-1}~~ {\rm and} ~~  C_{n,k} = C_{n,n-k-1}.
\]
\end{itemize}
\indent
\begin{itemize}
\item[(ii)] $n \equiv$ 2 or 3 ($\bmod$ 4). 
In this case, 
$A \in E_{\rm e}(n,k)$ if and only if
$A^{\ast} \in E_{\rm o}(n,n-k-1)$, and 
$A \in E_{\rm o}(n,k)$ if and only if
$A^{\ast} \in E_{\rm e}(n,n-k-1)$, so we have
\[
     B_{n,k} = C_{n,n-k-1}~~ {\rm and} ~~ C_{n,k} = B_{n,n-k-1}.
\]
\end{itemize}
\indent
The values of $B_{n,k}$ and $C_{n,k}$ for small $n$ 
are shown in the next two tables. The integers in their top rows represent the values of $k$.
In Section 4 a formula for calculating
these numbers will be supplied by means of $A_{n,k}$ and $D_{n,k}$. \\
\\
\begin{tabular}{l|rrrrrrrrrrr}
         \noalign{\hrule height 0.8pt}
        $B_{n,k}$ & 0 & 1 & 2 & 3 & 4 & 5 & 6 & 7 & 8 & ~~9  \\
       \hline
         $n = 2$ & 0 & 1 \\
         $n = 3$ & 0 & 2 & 1  \\
         $n = 4$ & 1 & 5 & 5 & 1 \\
         $n = 5$ & 1 & 14 & 30 & 14 & 1 \\
         $n = 6$ & 0 & 28 & 155 & 147 & 29 & 1 \\
         $n = 7$ & 0 & 56 & 605 & 1208 & 586 & 64 & 1 \\ 
         $n = 8$ & 1 & 127 & 2133 & 7819 & 7819 & 2133 & 127 & 1 \\ 
         $n = 9$ & 1 & 262 & 7288 & 44074 & 78190 & 44074 & 7288 & 262 &1 \\ 
         $n = 10$ & 0 & 496 & 23947 & 227623 & 655039 & 655315 & 227569 & 
         23893 & 517&~~1 \\ 
\noalign{\hrule height 0.8pt}
    \end{tabular}\\
\\
\\
\begin{tabular}{l|rrrrrrrrrr}
         \noalign{\hrule height 0.8pt}
        $C_{n,k}$ & 0 & 1 & 2 & 3 & 4 & 5 & 6 & 7 & 8 & ~~9  \\
       \hline
         $n = 2$ & 1 & 0 \\
         $n = 3$ & 1 & 2 & 0  \\
         $n = 4$ & 0 & 6 & 6 & 0 \\
         $n = 5$ & 0 & 12 & 36 & 12 & 0 \\
         $n = 6$ & 1 & 29 & 147 & 155 & 28 & 0 \\
         $n = 7$ & 1 & 64 & 586 & 1208 & 605 & 56 & 0 \\ 
         $n = 8$ & 0 & 120 & 2160 & 7800 & 7800 & 2160 & 120 & 0 \\ 
         $n = 9$ & 0 & 240 & 7320 & 44160 & 78000 & 44160 & 7320 & 240 & 0 \\ 
         $n = 10$ & 1 & 517 & 23893 & 227569 & 655315 & 655039 & 227623
         & 23947 & 496 & ~~0 \\ 
\noalign{\hrule height 0.8pt}
 \end{tabular} \\
\\
\\
\noindent
 {\bf 3.  The Case of Odd {\boldmath $n$}}\\
\\
\indent
Throughout this section we assume that $n$ is an odd integer. 
In this case, the orbit of a permutation of $E_{\rm e}^{-}(n,k)$ 
under $\sigma$ is entirely contained in $E_{\rm e}^{-}(n,k)$ and similarly 
for $E_{\rm e}^{+}(n,k)$, as was shown in Section 1. 
Here we mainly deal only with the set $E_{\rm e}(n,k)$
and its cardinality $B_{n,k}$, for the same arguments can also be applied to
$E_{\rm o}(n,k)$ and its cardinality $C_{n,k}$. \\
\indent
It was shown in [9] that to each permutation $A$ there corresponds a smallest positive 
integer $\pi(A)$ such that $\sigma^{\pi(A)}A = A$, which is called the 
{\it period} of $A$. Its trace 
\[
           \{\sigma A, \sigma^2 A, \ldots, \sigma^{\pi(A)}A=A \} 
\]
is called the {\it orbit} of $A$.
Also there it was shown that the period satisfies the relation
\begin{eqnarray}
     \pi(A) = 
        \left\{\begin{array}{rl}
                    (n-k)\gcd(n,\pi(A))  & \mbox{ if  $A \in E^{-}(n,k)$, } \\
                    (k+1)\gcd(n,\pi(A)) & \mbox{ if  $A \in E^{+}(n,k)$.} \\
                    \end{array}  \right.
\end{eqnarray}
\indent
It follows from (2) that the period of a permutation $A \in E(n,k)$ is 
either $d(n - k)$ or $d(k+1)$ for a positive divisor $d$ of $n$, i.e., $d = \gcd(n,\pi(A))$, 
although there may be no permutations having such periods for some divisors. 
In this paper, divisors of $n$ always mean positive divisors. \\
\indent
For a divisor $d$ of $n$, 
we denote by $\alpha_d^{k}$ the number of orbits of
period $d(n-k)$ in $E_{\rm e}^{-}(n,k)$  and 
by $\beta_d^{k}$ that of orbits of period $d(k+1)$ in $E_{\rm e}^{+}(n,k)$. 
In the case of odd $n$ the next theorem plays a fundamental role. \\
\\
\noindent
{\bf Theorem 3.1.}
Let $n$ be an odd integer and let $k$ be an integer satisfying 
$1 \le k \le n-1$.  Then it follows that
\begin{eqnarray}
         B_{n-1,k-1}  =  \sum_{d|n} d \alpha_d^{k},\\
         B_{n-1,k}  =   \sum_{d|n} d \beta_d^{k}, \\
         B_{n,k}  =  \sum_{d|n}d\{ (n-k) \alpha_d^{k} + (k+1)\beta_d^{k}\}. 
\end{eqnarray}
\\
\noindent
{\it Proof.} First let us consider permutations in $E_{\rm e}^{-}(n,k)$. Since
each orbit contains at least one canonical permutation, it suffices to deal only with
canonical ones in counting orbits.
If $A=1a_2a_3 \cdots a_n \in E_{\rm e}^{-}(n,k)$, we see that
$(a_2 -1)(a_3 -1) \cdots (a_n-1) \in E_{\rm e}(n-1,k-1)$, since
\[
{\rm inv}((a_2 -1)(a_3 -1) \cdots (a_n-1)) = {\rm inv}(A) 
\]
and one is deleted. Therefore,
there are $B_{n-1,k-1}$ canonical permutations in $E_{\rm e}^{-}(n,k)$. 
It follows from (2) that the period of a permutation $A \in E_{\rm e}^{-}(n,k)$ 
is equal to $d(n-k)$ for a divisor $d$ of $n$. 
There exist $n$ canonical permutations in
$\{ \sigma A, \sigma^2 A, \ldots, \sigma^{n(n-k)}A=A \}$ due to [9, Corollary 2], 
and hence each orbit $\{ \sigma A, \sigma^2 A, \ldots, \sigma^{d(n-k)}A=A \}$
of a permutation $A$ with period $d(n-k)$ contains exactly $d$ canonical permutations. 
This follows from the fact that the latter repeates itself $n/d$ 
times in the former. 
Since there exist $\alpha_d^{k}$ orbits of period $d(n-k)$ for each divisor $d$ of $n$,
classifying all canonical permutations of $E_{\rm e}^{-}(n,k)$ into orbits leads 
us to (3). \\
\indent
The proof of (4) is similar.  
To do this we consider permutations in $E_{\rm e}^{+}(n,k)$. 
If $A=a_2a_3 \cdots a_n1 \in E_{\rm e}^{+}(n,k)$, we see that
$(a_2 -1)(a_3 -1) \cdots (a_n-1) \in E_{\rm e}(n-1,k)$, since
\[
{\rm inv}((a_2 -1)(a_3 -1) \cdots (a_n-1)) = {\rm inv}(A) - (n-1)
\]
and $n-1$ is an even number by assumption. Therefore, the set of all canonical permutations 
in $E_{\rm e}^{+}(n,k)$ has cardinality $B_{n-1,k}$. 
Again using (2), the period of a permutation $A \in E_{\rm e}^{+}(n,k)$ is equal to $d(k+1)$ for a divisor $d$ of $n$. 
By [9, Corollary 2] there exist $n$ canonical permutations in
$\{ \sigma A, \sigma^2 A, \ldots, \sigma^{n(k+1)}A=A \}$ 
and hence, as above, there exist exactly $d$ such permutations in each orbit
$\{ \sigma A, \sigma^2 A, \ldots, \sigma^{d(k+1)}A=A \}$ of a permutation $A$ 
with period $d(k+1)$.
There exist $\beta_d^{k}$ orbits of period $d(k+1)$ 
for each divisor $d$ of $n$. Hence, we can obtain (4) by classifying all canonical
permutations in $E_{\rm e}^{+}(n,k)$ into orbits. \\
\indent
Considering the numbers of orbits and periods,
we see that the cardinalities of $E_{\rm e}^{\pm}(n,k)$ 
are obtained by 
\begin{eqnarray}
    |E_{\rm e}^{-}(n,k)|=  \sum_{d|n} d(n-k) \alpha_d^{k}
   ~~ {\rm and}~~ |E_{\rm e}^{+}(n,k)|= \sum_{d|n} d (k+1) \beta_d^{k}.
\end{eqnarray}  
Since the set $E_{\rm e}(n,k)$ is a disjoint union 
of $E_{\rm e}^{-}(n,k)$ and $E_{\rm e}^{+}(n,k)$, we conclude that
\begin{eqnarray*}
 B_{n,k} = |E_{\rm e}^{-}(n,k)|+|E_{\rm e}^{+}(n,k)|
     = \sum_{d|n} d(n-k) \alpha_d^{k} +\sum_{d|n} d (k+1) \beta_d^{k},
\end{eqnarray*}
which proves (5). \\
\\
\indent
Let us denote by $\gamma_d^{k}$ the number of orbits of
period $d(n-k)$ in $E_{\rm o}^{-}(n,k)$  and 
by $\delta_d^{k}$  that of orbits of period $d(k+1)$ in $E_{\rm o}^{+}(n,k)$. 
When $n$ is odd, analogous relations to (3)-(6) hold for $C_{n,k}$, 
$\gamma_d^{k}$ and $\delta_d^{k}$, since the orbit of a permutation of 
$E_{\rm o}^{\pm}(n,k)$ under $\sigma$ is also contained in 
$E_{\rm o}^{\pm}(n,k)$. We state them for the sake of completeness:
\begin{eqnarray*}
         C_{n-1,k-1}  =  \sum_{d|n} d \gamma_d^{k};\\
         C_{n-1,k}  =   \sum_{d|n} d \delta_d^{k}; \\
         C_{n,k}  =  \sum_{d|n}d\{ (n-k) \gamma_d^{k} + (k+1)\delta_d^{k}\}; 
\end{eqnarray*}  
and
\begin{eqnarray*}
    |E_{\rm o}^{-}(n,k)|=  \sum_{d|n} d(n-k) \gamma_d^{k}
   ~~ {\rm and}~~ |E_{\rm o}^{+}(n,k)|= \sum_{d|n} d (k+1) \delta_d^{k}.
\end{eqnarray*}
\indent
Making use of (3) and (4), we see that both cardinalities in (6) can be written simply
by $B_{n,k}$ and their counterparts for $E_{\rm o}(n,k)$ also follow from the above 
relations in a similar manner. \\
\\
\noindent
{\bf Corollary 3.2.}  When $n$ is odd, the cardinalities of $E^{\pm}_{\rm e}(n,k)$ 
and $E^{\pm}_{\rm o}(n,k)$
are given by
\begin{itemize}
\item[(i)] $|E^{-}_{\rm e}(n,k)|= (n-k)B_{n-1,k-1}$ ~~  and ~~ 
          $|E^{-}_{\rm o}(n,k)|= (n-k)C_{n-1,k-1}$ 
            ~ {\rm($1 \le k \le n-1$),} 
\item[(ii)] $|E^{+}_{\rm e}(n,k)| = (k+1)B_{n-1,k}$~~ and ~~ 
            $|E^{+}_{\rm o}(n,k)| = (k+1)C_{n-1,k}$ 
            ~ {\rm ($0 \le k \le n-2$). }\\
\end{itemize}

From these equalities we can obtain the following two corollaries. The relations in Corollary 3.3 have the same form as the recurrence relation for classical Eulerian numbers $A_{n,k}$; 
\begin{eqnarray}
A_{n,k} =  (n-k)A_{n-1,k-1} + (k+1)A_{n-1,k}.
\end{eqnarray}
The formula for $C_{n,k}$ can also be obtained from that for $B_{n,k}$ 
using the equality
$A_{n,k} = B_{n,k} + C_{n,k}$. \\

\noindent
{\bf Corollary 3.3.} 
When $n$ is odd, the following relations hold for $B_{n,k}$ and 
$C_{n,k}$: 
\begin{eqnarray}
B_{n,k} =  (n-k)B_{n-1,k-1} + (k+1)B_{n-1,k}; \\
C_{n,k} = (n-k)C_{n-1,k-1} + (k+1)C_{n-1,k}.
\end{eqnarray}
\noindent
{\bf Corollary 3.4.}  
When $n$ is odd, the following relations hold: 
\begin{itemize}
\item[(i)] $|E^{-}_{\rm e}(n,k)|-|E^{-}_{\rm o}(n,k)|= (n-k)D_{n-1,k-1}$
        ~ ($1 \le k \le n-1$); 
\item[(ii)] $|E^{+}_{\rm e}(n,k)|-|E^{+}_{\rm o}(n,k)|= (k+1)D_{n-1,k}$
       ~ ($0 \le k \le n-2$). \\
\end{itemize}
\noindent
 {\bf 4.  Recurrence Relation for Signed Eulerian Numbers $D_{n,k}$}\\
\\
\indent
When $n$ is even, equality (8) nor (9) does not hold,
as is seen from the tables of Section 2. 
For example, an odd integer $C_{10,4}$
cannot be written as a linear sum of $C_{9, k}$'s or $B_{9, k}$'s $(1 \le k \le 7)$ with 
integral coefficients, since they are all even.
Therefore, in reality (8) nor (9) does not
provide a recurrence relation of the numbers $B_{n,k}$ or $C_{n,k}$. \\
\indent
As for the differences $D_{n,k}=B_{n,k}-C_{n,k}$, however, 
their recurrence relation was conjectured in [7] and an analytic proof for it 
was given in [1]. 
In our notation it is described as the next theorem, 
for which we provide another proof 
from a combinatorial point of view. Notice that there is a different flavor in the case
of even $n$.\\
\\
\noindent
{\bf Theorem 4.1.} 
The recurrence relation for $D_{n,k}$ is given by
\begin{eqnarray}
     D_{n,k} = 
        \left\{\begin{array}{cl}
                    (n-k)D_{n-1,k-1}+ (k+1)D_{n-1,k} & \mbox{if }~n~ \mbox{ is odd},\\
                    D_{n-1,k-1} - D_{n-1,k} &  \mbox{if }~ n ~\mbox{ is even}. \\
                    \end{array}  \right.
\end{eqnarray}
\\
\noindent
{\it Proof.} 
The first part of this relation follows immediately from (8) and (9) of Corollary 3.3. 
Assuming that $n$ is even, we show the second part by means of
the operator $\sigma$. \\
\indent
Recall that when $n$ is even, the operator $\sigma$ may change the parity of permutations of $E(n,k)$ 
and it is a bijection on $E_{\rm e}^{-}(n,k) \cup E_{\rm o}^{-}(n,k)$ and on
$E_{\rm e}^{+}(n,k) \cup E_{\rm o}^{+}(n,k)$. \\
\indent
First let us consider permutations $A=a_1a_2 \cdots a_n$ in $E_{\rm e}^{-}(n,k) \cup E_{\rm o}^{-}(n,k)$ and divide all
permutations in $E_{\rm e}^{-}(n,k) \cup E_{\rm o}^{-}(n,k)$ into the following two types:
\begin{itemize}
\item[(i)]  $A=a_1a_2 \cdots a_{n-1}n$, where $a_1a_2 \cdots a_{n-1}$ is a permutation of $[n-1]$;
\item[(ii)]  $A=a_1a_2 \cdots a_n$ with $a_1 < a_n$, where $a_i = n$ for some $i$ $(2 \le i \le n-1)$.
\end{itemize}
Suppose $A \in E_{\rm e}^{-}(n,k)$. If $A$ is of type (i), then $\sigma A$ remains an even permutation, 
since ${\rm inv}(\sigma A) = {\rm inv}(A)$.
We see that the cardinality of permutations of type (i) is $B_{n-1,k-1}$, since $A$ is even and $n$ is the last entry. 
However, if $A\in E_{\rm e}^{-}(n,k)$ is of type (ii), then we have $\sigma A \in E_{\rm o}^{-}(n,k)$ by (1).
Therefore, the cardinality of permutations of type (ii) in $E_{\rm e}^{-}(n,k)$ is 
\[
             |E_{\rm e}^{-}(n,k)| - B_{n-1,k-1},
\]
and precisely so many permutations change the parity from even to odd under $\sigma$.\\
\indent
Simliarly, suppose $A \in E_{\rm o}^{-}(n,k)$. If $A$ is of type (i), then $\sigma A$ 
remains an odd permutation.
We see that the cardinality of permutations of type (i)  is $C_{n-1,k-1}$.
If $A \in E_{\rm o}^{-}(n,k)$ is of type (ii) by (1), then we have $\sigma A \in E_{\rm e}^{-}(n,k)$. 
The cardinality of permutations of type (ii) in $E_{\rm o}^{-}(n,k)$ is 
\[
             |E_{\rm o}^{-}(n,k)| - C_{n-1,k-1},
\]
and precisely so many permutations change the parity from odd to even under $\sigma$.\\
\indent
Since $\sigma$ is a bijection on $E_{\rm e}^{-}(n,k) \cup E_{\rm o}^{-}(n,k)$, both cardinalities must
be equal. Hence we obtain 
\begin{eqnarray}
|E_{\rm e}^{-}(n,k)| - |E_{\rm o}^{-}(n,k)| = B_{n-1,k-1} - C_{n-1,k-1} = D_{n-1,k-1}.
\end{eqnarray}
\indent
Next let us consider permutations $A=a_1a_2 \cdots a_n$ in $E_{\rm e}^{+}(n,k) \cup E_{\rm o}^{+}(n,k)$ and divide all
permutations in $E_{\rm e}^{+}(n,k) \cup E_{\rm o}^{+}(n,k)$ into the following two types:
\begin{itemize}
\item[(iii)]  $A=na_1a_2 \cdots a_{n-1}$, where $a_1a_2 \cdots a_{n-1}$ is a permutation of $[n-1]$;
\item[(iv)]  $A=a_1a_2 \cdots a_n$ with $a_1 > a_n$, where $a_i = n$ for some $i$ $(2 \le i \le n-1)$.
\end{itemize}
If $A \in E_{\rm e}^{+}(n,k)$ is of type (iii), then $\sigma A$ remains an even permutation.
We see that the cardinality of permutations of type (iii) is $C_{n-1,k}$, 
since ${\rm inv}(A) -{\rm inv}(a_1a_2 \cdots a_{n-1})= n-1$ and $n-1$ is odd. 
However, if $A \in E_{\rm e}^{+}(n,k)$ is of type (iv), then we have $\sigma A \in E_{\rm o}^{+}(n,k)$ by (1).
The cardinality of permutations of type (iv) in $E_{\rm e}^{+}(n,k)$ is 
\[
             |E_{\rm e}^{+}(n,k)| - C_{n-1,k},
\]
and precisely so many permutations change the parity from even to odd under $\sigma$.\\
\indent
Similarly, if $A \in E_{\rm o}^{+}(n,k)$ is of type (iii), then $\sigma A$ remains an odd permutation.
We see that the cardinality of permutations of type (iii) is $B_{n-1,k}$ as above.
If $A \in E_{\rm o}^{+}(n,k)$ is of type (iv), then we have $\sigma A \in E_{\rm e}^{+}(n,k)$.
The cardinality of permutations of type (iv) in $E_{\rm o}^{+}(n,k)$ is 
\[
             |E_{\rm o}^{+}(n,k)| - B_{n-1,k},
\]
and precisely so many permutations change the parity from odd to even under $\sigma$.\\
\indent
Since $\sigma$ is a bijection on $E_{\rm e}^{+}(n,k) \cup E_{\rm o}^{+}(n,k)$, both cardinalities       
must be equal. Hence we obtain 
\begin{eqnarray}
|E_{\rm e}^{+}(n,k)| - |E_{\rm o}^{+}(n,k)| = -B_{n-1,k} + C_{n-1,k} = -D_{n-1,k}.
\end{eqnarray}
Adding (11) and (12) yields 
\[
     B_{n,k}- C_{n,k} =  D_{n,k} = D_{n-1,k-1}-D_{n-1,k},
\]
which is the required relation. This completes the proof. \\
\\
\indent
Symmetry properties for them follow from the relations presented in Section 2:
\indent
\begin{itemize}
\item[(i)] For $n \equiv$ 0  or 1 ($\bmod$ 4),  $D_{n,k} = D_{n,n-k-1}$;
\indent
\item[(ii)] For $n \equiv$ 2 or 3 ($\bmod$ 4), $D_{n,k} = - D_{n,n-k-1}$.
\end{itemize}
A table for the values of $D_{n,k}$ is given below. \\
\\
\begin{tabular}{l|rrrrrrrrrrr}
         \noalign{\hrule height 0.8pt}
        $D_{n,k}$ & 0 & 1 & 2 & 3 & 4 & 5 & 6 & 7 & ~8 & ~~9  \\
       \hline
         $n = 2$ & $-1$ & 1 \\
         $n = 3$ & $-1$ & 0 & 1  \\
         $n = 4$ & 1 & $-1$ & $-1$ & 1 \\
         $n = 5$ & 1 & 2 & $-6$ & 2 & 1 \\
         $n = 6$ & $-1$ & $-1$ & 8 & $-8$ & 1 & 1 \\
         $n = 7$ & $-1$ & $-8$ & 19 & 0 & $-19$ & 8 & 1 \\ 
         $n = 8$ & 1 & 7 & $-27$ & 19 & 19 & $-27$ & 7 & 1 \\ 
         $n = 9$ & 1 & 22 & $-32$ & $-86$ & 190 & $-86$ & $-32$ & 22 &1 \\ 
         $n = 10$ & $-1$ & $-21$ & 54 & 54 & $-276$ & 276 & $-54$ & 
         $-54$ & ~21&~~1\\ 
\noalign{\hrule height 0.8pt}
    \end{tabular}\\
\\
\indent
Thus the values of $B_{n,k}$ and $C_{n,k}$ can be known through
\[
     B_{n,k}= \frac{A_{n,k}+D_{n,k}}2, ~~   C_{n,k}= \frac{A_{n,k}-D_{n,k}}2, 
\]
using $A_{n,k}$ and $D_{n,k}$ that are calculated according to 
the respective recurrence relations (7) and (10). 
From these equalities, we can obtain the expressions of $B_{n,k}$ and $C_{n,k}$ by means of
$B_{n-1,k}$'s and $C_{n-1,k}$'s in the case of even $n$, which is a counterpart of Corollary 3.3.\\
\\
\noindent
{\bf Corollary 4.2.} 
When $n$ is even, the following relations hold for $B_{n,k}$ and 
$C_{n,k}${\rm :} 
\begin{eqnarray*}
2B_{n,k} =  (n-k+1)B_{n-1,k-1} + kB_{n-1,k}+(n-k-1)C_{n-1,k-1} + (k+2)C_{n-1,k}; \\
2C_{n,k} =  (n-k+1)C_{n-1,k-1} + kC_{n-1,k}+(n-k-1)B_{n-1,k-1} + (k+2)B_{n-1,k} .
\end{eqnarray*} 
\\
\indent
From (11) and (12) we get a counterpart of Corollary 3.4. \\
\\
\noindent
{\bf Corollary 4.3.} 
When $n$ is even, the following relations hold:
\begin{itemize}
\item[(i)] $|E^{-}_{\rm e}(n,k)|-|E^{-}_{\rm o}(n,k)|= D_{n-1,k-1}$
        ~ ($1 \le k \le n-1$); 
\item[(ii)] $|E^{+}_{\rm e}(n,k)|-|E^{+}_{\rm o}(n,k)|= -D_{n-1,k}$
       ~ ($0 \le k \le n-2$). \\
\end{itemize}
\noindent
{\bf 5.  Numbers of Orbits and Applications}\\
\\
\indent
Again $n$ is assumed to be an odd integer. In this section we derive the numbers of 
orbits of particular types, and moreover deduce divisibility properties for $B_{n,k}$, 
$C_{n,k}$, $D_{n,k}$ and related numbers by prime powers from them. \\
\indent
For a positive integer $\ell$ with $\gcd(\ell, n) = 1$, a canonical permutation 
of $[n]$ of the form
\[
P_n^{\ell} = 1(1+\ell)(1+2\ell) \cdots (1+(n-1)\ell)
\]
 can be defined, where $\ell, 2\ell, \ldots, (n-1)\ell$ represent numbers modulo $n$. 
According to whether $P_n^{\ell}$ is an even or odd permutation, let us put
\begin{eqnarray*}
     \epsilon^{\ell}_n = 
        \left\{\begin{array}{rl}
                    1 & \mbox{~if ~$P_n^{\ell}$ is even, } \\
                    0  & \mbox{~if ~$P_n^{\ell}$ is odd.} \\
                    \end{array}  \right.
\end{eqnarray*}
\noindent
{\bf Theorem 5.1.} 
Let $n$ be an odd integer and let $k$ be an integer such that 
$1 \le k \le n-1$. 
\begin{itemize}
\item[(i)]  If a divisor $d$ of $n$ satisfies $\gcd(k, n/d) > 1,$ 
    then $\alpha_d^k =\gamma_d^k = 0$. 
\item[(ii)] If $\gcd(k, n) = 1,$ then $\alpha_1^{k} =  \epsilon^{n-k}_n$
and $\gamma_1^{k} = 1- \epsilon^{n-k}_n$.
\end{itemize}
\noindent
{\it Proof.} 
In order to prove (i), suppose $A$ is a permutation that belongs to 
$E_{\rm e}^{-}(n,k)$. From (2) its period $\pi(A)$ satisfies
$\pi(A) = (n-k)\gcd(n, \pi(A))$. Then, putting 
$d = \gcd(n, \pi(A))$, we have $\pi(A) = d(n-k)$ and $d = \gcd(n, d(n-k))$, 
which implies $\gcd(n-k, n/d) = 1$ or $\gcd(k, n/d) = 1$. 
Consequently, we see that there exist no permutations of 
period $d(n-k)$, i.e., $\alpha_d^k = 0$, if a divisor $d$ of $n$ satisfies 
$\gcd(k, n/d) > 1$. 
The same arguments can be applied to permutations in $E_{\rm o}^{-}(n,k)$ 
and we obtain the assertion that
$\gamma_d^k = 0$ if $d$ satisfies $\gcd(k, n/d) > 1$.  \\
\indent
Next suppose $\gcd(k, n) = 1$. By [9, Theorem 7] we see that there exists a unique 
orbit of period $n-k$ in $E_{\rm e}^{-}(n,k) \cup E_{\rm o}^{-}(n,k)$, 
which contains only one canonical permutation
$P_n^{n-k}$. Hence, if it is an even permutation, then we have
$\alpha_1^{k} = 1$ and $\gamma_1^{k} = 0$. Otherwise, 
$\alpha_1^{k} = 0$ and $\gamma_1^{k} = 1$.  This completes the proof. \\
\\
\indent
From Theorem 5.1 we can derive a criterion under which $B_{n -1, k-1}$ and
$C_{n-1,k-1}$ is divisible by a prime power. \\
\\
\noindent
{\bf Corollary 5.2.}
Let $p$ be a prime and let an odd integer $n$ be divisible by $p^m$ for a positive integer $m$. 
If $k$ is divisible by $p$, then $B_{n -1, k-1}$, $C_{n-1,k-1}$ and $D_{n-1,k-1}$ 
are also divisible by $p^m$.
\\
\\
\noindent
{\it Proof.}  
Without loss of generality we can assume 
that $m$ is the largest integer for which $p^m$ divides $n$. 
Suppose $k$ is a multiple of $p$.
In Theorem 5.1 we have seen that $\alpha_d^k = 0$ for a divisor $d$ of $n$ 
such that $\gcd(k, n/d)>1$.
On the other hand, it follows that a divisor $d$ for which $\gcd(k, n/d) =1$ 
must be a multiple of $p^m$, since $k$ is a multiple of $p$. Therefore, 
equality (3) of Theorem 3.1 implies that $B_{n -1, k -1}$ is divisible by $p^m$.
Similarly, $C_{n -1, k -1}$ is also divisible by $p^m$, if $k$ is a multiple of $p$.\\
\\
\indent
The final corollary easily follows from Corollaries 3.2 and 5.2. \\
\\
\noindent
{\bf Corollary 5.3.} 
Under the same assumptions as Corollary 5.2 it follows that:
\begin{itemize}
\item[ (i)] If $k$ is divisible by $p^i$ for some $i$ $(1 \le i \le m)$, 
then $|E^{-}_{\rm e}(n,k)|$ and $|E^{-}_{\rm o}(n,k)|$ are divisible by $p^{m+i}$.
\item[(ii)]  If $k+1$ is divisible by $p^i$ for some $i$ ($i \ge 1$),
then $|E^{+}_{\rm e}(n,k)|$ and $|E^{+}_{\rm o}(n,k)|$ are divisible by $p^{m+i}$.
\end{itemize}
\noindent
{\bf References}
\begin{itemize}
\item[1.]   J. D\'esarm\'enien and D. Foata, 
The signed Eulerian numbers, Discrete Math. 99 (1992) 49-58.
\item[2.]  D. Foata, and M.-P.  Sch\"utzenberger, 
Th\'eorie G\'eom\'etrique des Polyn\^omes Eul\'eriens,  Lecture Notes in
Mathematics, Vol. 138, Springer-Verlag, Berlin, 1970.
\item[3.] R. L. Graham, D. E. Knuth and  O. Patashnik,  Concrete Mathematics,
Addison-Wesley, Reading, 1989.
\item[4.]  A. Kerber, Algebraic Combinatorics Via Finite Group Actions. 
BI-Wissenschaftsverlag, Mannheim, 1991. 
\item[5.]   D. E.  Knuth, The Art of Computer Programming, Vol. 3, Sorting and
Searching. Addison-Wesley, Reading, 1973.
\item[6.]   L. Lesieur and J.-L. Nicolas,  On the Eulerian numbers
    $M_n =\max_{1\le k \le n}A(n,k)$,  European J. Combin. 13 (1992)  379-399.
\item[7.]   J.-L. Loday,  Op\'erations sur l'homologie cyclique des alg\`ebres
commutatives, Invent. Math. 96 (1989) 205-230. 
\item[8.]   R. Mantaci,  Binomial coefficients and anti-excedances of even 
permutations: A combinatorial proof, J. of Comb. Theory (A)  63
(1993) 330-337.
 \item[9.]    S. Tanimoto,  An operator on permutations and its 
application to Eulerian numbers, European J. Combin. 22 (2001)  569-576.
 \item[10.]    S. Tanimoto,  A study of Eulerian numbers by means of an operator on permutations, 
European J. Combin. 24 (2003) 34-44.
 \item[11.]    S. Tanimoto, On the numbers of orbits of permutations under an operator related to Eulerian numbers, 
Annals of Combin. 8 (2004) 239-250.
\end{itemize}

\end{document}